\documentclass[11pt]{article}
\usepackage{amsfonts}
\usepackage{amsmath}
\usepackage{mathrsfs}
\usepackage{color}
\usepackage{mathrsfs,amscd,amssymb,amsthm,amsmath,bm,graphicx,psfrag,subfigure,url}

\makeatletter

\usepackage{indentfirst}

\def \[{\begin{equation}}
\def \]{\end{equation}}

\newtheorem{thm}{Theorem}[section]

\newtheorem{claim}{Claim}
\newtheorem{lem}[thm]{Lemma}
\newtheorem{cor}[thm]{Corollary}

\newenvironment{wst}
{\setlength{\leftmargini}{1.5\parindent}
 \begin{itemize}
 \setlength{\itemsep}{-1.1mm}}
{\end{itemize}}

\begin{document}

\title{\bf Further results on the expected hitting time, the cover cost and the related invariants of graphs \footnote{S. L. acknowledges the financial support from the National Natural Science Foundation of China (Grant Nos. 11671164, 11271149).}}

\author{Jing Huang $^a$, \ Shuchao Li $^{a,}$\footnote{Corresponding author.  }\, ,\ Zheng Xie $^b$}

\date{}

\maketitle

\begin{center}
 Faculty of Mathematics and Statistics, Central China Normal
University, \\ Wuhan 430079, P.R. China\\[5pt]
College of Science, National University of Defense Technology, \\
Changsha 410073, P.R. China\\[5pt]
jhuangmath@sina.com (J.~Huang), lscmath@mail.ccnu.edu.cn (S.C.~Li), xiezheng81@nudt.edu.cn (Z.~Xie)

\end{center}

%%%%%%%%%%%%%

\begin{abstract}
A close relation between hitting times of the
simple random walk on a graph, the Kirchhoff index, resistance-centrality, and related invariants of unicyclic graphs is displayed. Combining with the graph transformations and some other techniques, sharp upper and lower bounds on the cover cost (resp. reverse cover cost) of a vertex in an $n$-vertex unicyclic graph are determined. All the corresponding extremal graphs are identified, respectively.

%Let $\mathcal{U}_n$ be the set of all unicyclic graphs on $n$ vertices and $\mathscr{U}_{n,l}$ be the set of all unicyclic graphs on $n$ vertices containing the unique cycle $C_l$. The \textit{cover cost} (resp. \textit{reverse cover cost}) of a vertex $x$ in $G$ is defined as $CC(x)=\sum_{y\in V_G}H_{xy}$ (resp. $RC(x)=\sum_{y\in V_G}H_{yx}$ ), where the expected hitting time $H_{xy}$ is the expected number of steps it takes a simple random walk on a graph $G$ to go from $x$ to $y$. The Kirchhoff index of a graph $G$ is denoted by $K\!f(G)=\sum_{\{x,y\}\subseteq V_G}r(x,y)$, where $r(x,y)$ is the resistance distance between vertices $x$ and $y$ in $G$. In this paper, we first exhibit some close connections between the expected hitting time, the Kirchhoff index and some other related graph invariants of unicyclic graphs. Then we establish the sharp lower and upper bounds on $CC(x)$ of a graph $G$ in $\mathscr{U}_{n,l}$ (resp. $\mathcal{U}_n$) and characterize the corresponding extremal graphs, respectively. Finally,  the sharp lower and upper bounds on $RC(x)$ of graph $G$ in $\mathscr{U}_{n,l}$ (resp. $\mathcal{U}_n$) are given, and the corresponding extremal graphs are determined, respectively.
\end{abstract}

\vspace{2mm} \noindent{\bf Keywords}: Random walk; Hitting time;  Kirchhoff index; Cover cost; Reverse cover cost

\vspace{2mm}


\setcounter{section}{0}
\section{Introduction}\setcounter{equation}{0}

We will start with introducing some background information that will lead to our main results. Some important previously established facts will also be presented.

\subsection{Background}
All graphs considered in this paper are simple and undirected. Let $G=(V_G, E_G)$ be a graph with $V_G$ the vertex set and $E_G$ the edge set.  We call $n=|V_G|$ the \textit{order} of $G$ and $m=|E_G|$ the \textit{size} of $G$. The \textit{neighborhood} of a vertex $x$, written by $N(x)$, is the set of vertices adjacent to $x$ in $G$. The degree of $x$ is $d(x)=|N(x)|$. A vertex of a graph $G$ is called a \textit{pendant vertex} if it is of degree $1$. The \textit{distance} between vertices $x$ and $y$, denoted by $d(x,y)$, is the length of a shortest path connecting them and the \textit{eccentricity} $\varepsilon(x)$ of a vertex $x$ is the distance between $x$ and a furthest vertex from $x$ in $G$. We call a simple graph $G$ a \textit{unicyclic graph} if it is connected with $|V_G|=|E_G|$.  We follow the notation and terminologies in \cite{D-I} except if otherwise stated.

An important parameter called \textit{Wiener index}, $W(G),$ was defined by $W(G)=\sum_{\{x,y\}\subseteq V_G}d(x,y)$ in \cite{w1}. It has been extensively studied and has found applications in chemistry, communication theory, and elsewhere. The \textit{centrality} (also known as the \textit{transmission}) of a vertex $x$ in $G$ is defined as $D(x)=\sum_{y\in V_G}d(x,y),$ whereas the \textit{weighted centrality} of a vertex $x$ is defined as $D^w(x)=\sum_{y\in V_G}d(y)d(x,y)$. Then it is obvious that $W(G)=\frac{1}{2}\sum_{x\in V_G}D(x)$.

The \textit{adjacency matrix} $A(G)$ of $G$ is an $n\times n$ matrix with the $(x,y)$-entry equals to 1 if vertices $x$ and $y$ are adjacent and 0 otherwise. Let $D(G)={\rm diag}(d_1,d_2,\ldots,  d_n)$ be the diagonal matrix of vertex degrees. The (\textit{combinatorial}) \textit{Laplacian matrix} of $G$ is defined as $L(G)=D(G)-A(G).$

Based on the electrical network theory, Klein and Randi\'{c} \cite{13} proposed a new distance-based parameter, i.e., the \textit{resistance distance}, on a graph. The resistance distance between vertices $x$ and $y$, written by $r(x,y),$ is the effective resistance between them when one puts one unit resistor on every edge of a graph $G$. It is known that $r(x,y)\leq d(x,y)$ with equality if and only if $G$ is a tree. One famous resistance distance-based invariant called the \textit{Kirchhoff index}, $K\!f(G)$, was given by $K\!f(G)=\sum_{\{x,y\}\subseteq V_G}r(x,y)$ (see \cite{13}). This structure-descriptor can be expressed alternatively as
$$
K\!f(G)=\sum_{\{x,y\}\subseteq V_G}r(x,y)=n\sum_{i=2}^n\frac{1}{\mu_i},
$$
where $0=\mu_1 <\mu_2\leqslant \cdots \leqslant\mu_n$ are the eigenvalues of $L(G).$

As an analogue of the Kirchhoff index of $G$, Chen and Zhang \cite{3} proposed a novel resistance distance-based graph invariant, defined by $K\!f^*(G)=\sum_{\{x,y\}\subseteq V_G}d(x)d(y)r(x,y)$, which is called the \textit{multiplicative degree-Kirchhoff index} (see \cite{20}). Just as the relationship between the Kirchhoff index and the Laplacian spectrum, the multiplicative degree Kirchhoff index is closely related to the spectrum of the normalized Laplacian matrix $\mathcal{L}(G)$, which is defined as $\mathcal{L}(G)={D(G)}^{-\frac{1}{2}}L(G){D(G)}^{-\frac{1}{2}}.$ It can be expressed alternatively as
\begin{eqnarray}\label{eq:1.1a}
K\!f^*(G)=\sum_{\{x,y\}\subseteq V_G}d(x)d(y)r(x,y)=2m\sum_{i=2}^n\frac{1}{\lambda_i},
\end{eqnarray}
where $0=\lambda_1 <\lambda_2\leqslant \cdots \leqslant\lambda_n$ are the eigenvalues of $\mathcal{L}(G).$  More recently, another resistance distance-based graph invariant, namely the \textit{additive degree-Kirchhoff index} has been put forward in \cite{E2}. It is defined as
\begin{eqnarray}\label{eq:1.1b}
K\!f^+(G)=\sum_{\{x,y\}\subseteq V_G}\left(d(x)+d(y)\right)r(x,y).
\end{eqnarray}
There is extensive literature available on works related to $K\!f(G),\, K\!f^*(G)$ and $K\!f^+(G)$, one may be referred to \cite{a1,a3,a2,a5,a6} for more detailed information.

For a graph $G$, define the random walks on $G$ as the Markov chain $X_k, k\geq0,$ that from its current vertex $x$ jumps to its adjacent vertex with probability $1/d(x)$. The \textit{hitting time} (also known as the \textit{first passage time}) $T_{y}$ of the vertex $y$ is the
minimum number of jumps the walk needs to reach $y$, that is
$$
T_y=\inf\{k\geq0:X_k=y\}.
$$
The expected value of $T_{y}$ when the walk is started at the vertex $x$ is denoted by $H_{xy}$.

The hitting time of random walks is an important parameter of graphs \cite{v1,v2} and it has been studied extensively. The connections of eigenvalues and  eigenvectors with hitting time were studied in \cite{A3}. Relationships between hitting time and electrical networks were considered in \cite{v5,v4,v6,v7}. In 2013, Xu and Yau \cite{40} proposed the $R$-invariant and $Z$-invariant (also called Chung-Yau invariants in \cite{8}) and provided an explicit formula of hitting time in terms of Chung-Yau invariants and the number of spanning trees \cite{c1}. In 2016, Patel et al. \cite{36} provided a novel method for calculating the hitting time for a single random walker as well as the first analytic expression for calculating the hitting
time for multiple random walkers, which they denoted as the \textit{group hitting time}. A closed form solution for calculating the hitting time between specified nodes for both the single and multiple random walker cases were also presented. Recently, Chang and Xu \cite{8} used the Chung-Yau graph invariants to derive new explicit formulas and estimated for hitting times of random walks. And they also applied these invariants to study graphs with symmetric hitting times. We refer the readers to the nice survey \cite{A3} and \cite{9,b3,17,24,13,28,33} for more background of random walks and hitting times on graphs.

The \textit{cover cost} (see \cite{g1}) of a vertex $x$ in $G$ is defined as
$$
   CC(x)=\sum_{y\in V_G}H_{xy}.
$$
It is closely related to the \textit{cover time} of a graph, which is defined as the expected time for a random walk starting at $x$ to visit all vertices. There is a rather beautiful relationship between the cover cost and the Wiener index of a tree, which is expressed as $CC(x)+D(x)=2W(T)$ for every tree $T$ and every vertex $x\in V_{T}$. As an analogue of the cover cost, the  \textit{reverse cover cost} of a vertex $x$ in $G$, which is defined as
$$
  RC(x)=\sum_{y\in V_G}H_{yx}.
$$
It was proposed by Georgakopoulos and Wagner \cite{w2} in which they showed that $RC(x)+(2n-1)CC(x)=4(n-1)W(T)$ for every $n$-vertex tree $T$ with $x\in V_T.$ As well, they determined the extremal values of the hitting time, the cover cost, and the reverse cover cost for trees of given order. All the the corresponding  extremal graphs were characterized.

Apart from all the trees, any unicyclic graph is a connected graph with as small size as possible. Motivated from \cite{w2}, it is natural and interesting for us to consider the problems as above for unicyclic graphs. Our methods and technique are novel, which are completely different from those in \cite{w2}.

\subsection{Main results}

In this subsection we give necessary definitions and state the main results of the paper. If $x\in V_G,$ then $G-x$ denotes the graph obtained from $G$ by deleting the vertex $x$ and all its incident edges. If $xy\notin E(G),$ then $G+xy$ is a graph obtained from $G$ by adding an edge $xy.$ For $X\subseteq E_G$, $G-X$ denotes the graph obtained from $G$ by deleting all the edges in $X$. In particular, if $X=\{xy\},$  we write $G-xy$ for $G-X$. Denote by $C_n, P_n$ and $S_n$  the cycle, the path and the star of order $n$, respectively.

Recall that a unicyclic graph $G$ is a simple connected graph with $|V_G|=|E_G|$. For convenience, we may use the following notation in the whole context to represent a unicyclic graph: Let $G=U(C_l; T_1, T_2,\ldots,T_l)$ be an $n$-vertex unicyclic graph, where $C_l=v_1 v_2\ldots v_lv_1$ is the unique cycle contained in $G$ and $T_i$ is the component of $G-E_{C_l}$ containing $v_i.$ Denote $n_i=|V_{T_i}|$ for $1\leq i\leq l$. It is easy to see that $T_i$ is a tree for $1\leq i\leq l$ and $\sum_{i=1}^ln_i=n.$  We say that $T_i$ is \textit{trivial} if $n_i=1$.
\begin{figure}[h!]
\begin{center}
\psfrag{a}{$C_l$}\psfrag{b}{$n-l$}
\psfrag{c}{$S_n^l$}\psfrag{d}{$P_n^l$}
\includegraphics[width=90mm]{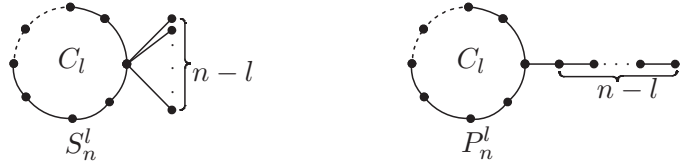} \\
\caption{Graphs $S_n^l$ and $P_n^l$} %\label{*}
\end{center}
\end{figure}

Let $\mathscr{U}_n$ be the set of all unicyclic graphs on $n$ vertices and $\mathscr{U}_{n,l}$ be the set of all $n$-vertex unicyclic graphs each of which contains the unique cycle $C_l$. Let $S_n^l$ denote the graph obtained  by attaching $n-l$ pendent vertices to exactly one vertex of $C_l$, and $P_n^l$ be the graph obtained by identifying one end vertex of $P_{n-l+1}$ with one vertex of $C_l.$ Graphs $S_n^l$ and $P_n^l$ are depicted in Fig. 1. It is obvious that $S_n^n=P_n^n=C_n.$

Taking advantage of the theory of the relationship between random walks and electrical networks \cite{010,v7}, it is natural to define the  \textit{resistance-centrality} and \textit{weighted resistance-centrality} (see \cite{w2}), respectively, as
$$
   R(x)=\sum_{y\in V_G}r(x,y), \ \ \ \ \ \ \ \ \ R^w(x)=\sum_{y\in V_G}d(y)r(x,y),
$$
which can be seen as the generalizations of $D(x)$ and $D^w(x)$. Obviously, $K\!f(G)=\frac{1}{2}\sum_{x\in V_G}R(x).$

Following the above notations our first main results read as follows.
\begin{thm}\label{th2}
Let $G=(C_l;T_1,T_2,\ldots, T_l)$ be in $\mathscr{U}_{n,l}$ with $x\in V_{T_k}.$ Then we have
\begin{wst}
\item[{\rm (i)}]  $CC(x)=2K\!f(G)+\sum_{i=1}^lD_{T_i}(v_i)-nd(x, v_k);$
\item[{\rm (ii)}] $RC(x)=2nR(x)-2K\!f(G)-\sum_{i=1}^lD_{T_i}(v_i)+nd(x, v_k).$
\end{wst}
\end{thm}
The following result is an immediately consequence of Theorem \ref{th2} (i).
\begin{cor}
Let $G=(C_l;T_1,T_2,\ldots, T_l)$ be in $\mathscr{U}_{n,l}$ with $x\in V_{T_i}, y\in V_{T_j}.$ Then $d(x,v_i)\leq d(y,v_j)$ if and only if $CC(x)\geq CC(y).$
\end{cor}

Our second main results in the following exhibit a explicit formula of the additive degree-Kirchhoff index (resp. multiplicative degree-Kirchhoff index) in terms of the Kirchhoff index, the order and the related invariants of unicyclic graphs.
\begin{thm}\label{th1}
Let $G=(C_l;T_1,T_2,\ldots, T_l)$ be in $\mathscr{U}_{n,l}.$ Then
\begin{wst}
\item[{\rm (i)}]  $K\!f^+(G)=4K\!f(G)+2\sum_{i=1}^lD_{T_i}(v_i)-n(n-l);$
\item[{\rm (ii)}] $K\!f^*(G)=4K\!f(G)+4\sum_{i=1}^lD_{T_i}(v_i)-(2n+1)(n-l).$
\end{wst}
\end{thm}
\noindent{\bf Remark 1.}\ \ Gutman et al. \cite{20,E2} characterized $n$-vertex unicyclic graphs having the minimum, and the second-minimum multiplicative degree Kirchhoff indices and additive degree-Kirchhoff indices, respectively. By Theorem \ref{th1}, we may give a new and much simpler proof for the main results in \cite{20,E2}, which are left to the readers.

Let
\[\label{eq:1.3}
f_1(n)=\left\{
           \begin{array}{llll}
            \frac{n^3-n^2+4n-6}{6}, & \hbox{if $4\leq n\leq 8;$} \\
            \frac{n^3-n}{6}, & \hbox{if $n=3,9,10;$}\\
             2n^2-5n-6, & \hbox{if $11\leq n\leq15;$}\\
             2n^2-\frac{16}{3}n-1, & \hbox{if $n\geq16,$}
             \end{array}
           \right.
\]
and
\[\label{eq:1.4}
  f_2(n)=\frac{2n^3+3n^2-37n+54}{6}.
\]
Then our third main results read as follows.
\begin{thm}\label{th4.50}
Let $G\in\mathscr{U}_n$ with $x\in V_G$. Then
$$
f_1(n)\leqslant CC(x)\leqslant f_2(n),
$$
where $f_1(n)$ and $f_2(n)$ are defined in $(\ref{eq:1.3})$ and $(\ref{eq:1.4})$, respectively. The right equality holds if and only if $G\cong P_n^3$ and $x$ lies on the unique cycle $C_3$. The left equality holds if and only if
\begin{wst}
\item[{\rm (i)}]for $4\leqslant n \leqslant 8,$  $G\cong S_n^{n-1}$ and $x$ is a pendant vertex;
\item[{\rm (ii)}]for $n=3,9,10,$  $G\cong  C_n;$
\item[{\rm (iii)}]for $11\leqslant n \leqslant 15,$  $G\cong S_n^4$ and $x$ is a pendant vertex;
\item[{\rm (iv)}]for $n\geqslant 16,$  $G\cong S_n^3$ with $x$ being a pendant vertex.
\end{wst}
\end{thm}

For our last main results, we concentrate on the reverse cover cost of unicyclic graphs. More precisely, sharp upper and lower bounds on $RC(x)$ of graph $G$ in $\mathscr{U}_n$ are established and the corresponding extremal graphs are determined, respectively.

For $X\subseteq V_G$ and $x\in V_G\setminus X$, define the \textit{distance} between $x$ and $X$ as $d(x,X)=\min\{d_G(x,y):\, y\in X\}.$
\begin{thm}\label{thm5}
Let $G\in\mathscr{U}_n$ with $x\in V_G$. Then
$$
  n+1\le RC(x)\leq \frac{n(n-1)(4n+1)}{6}-9.
$$
The left equality holds if and only if $G\cong S_n^3$ and $d(x)=n-1$, whereas the right equality holds if and only if $G\cong P_n^3$ and $x$ is the vertex such that $d(x,V_{C_3})=n-3.$
\end{thm}
In the rest of this section we recall some important known results. In Section~\ref{sec:2} we establish some technical lemmas that help us prove the main results. We present the proofs of Theorems \ref{th2} and \ref{th1} in Section~\ref{sec:3}. In Section~\ref{sec:4}, we give the proof for Theorem \ref{th4.50}, whereas in Section~\ref{sec:5}, we give the proof for Theorem \ref{thm5}.
\subsection{Preliminaries}

For the rest of our introduction we recall the following important facts.

\begin{lem}[\cite{v7}]\label{lem3.1}
Let $G$ be a simple connected graph on $m$ edges with $x,y\in V_G.$ Then
$$
H_{xy}=\frac{1}{2}\sum_{z\in V_G}d(z)(r(x,y)+r(z,y)-r(z,x))=mr(x,y)+\frac{1}{2}(R^w(y)-R^w(x)).
$$
\end{lem}
%The expected commute time between vertices $x$ and $y$ is defined by $H_{xy}(G)+H_{yx}(G).$ In view of \cite{A-K-C,13}, we know that there is an elegant relation between commute times and resistance distances, which can immediately follow from Lemma \ref{lem3.1}.
%\begin{lem}[\cite{A-K-C,13}]\label{lem3.2}
%Let $G$ be a connected graph $m$ edges with $x,y\in V_G.$ Then
%$$
%H_{xy}+H_{yx}=2mr(x,y).
%$$

%\end{lem}

\begin{lem}[\cite{w1}]\label{lem3.3}
Let $z$ be a cut vertex of a connected graph $G$ and $x, y$ be two vertices occurring in different components of $G-z$. Then
$$
r(x,y)=r(x,z)+r(z, y).
$$
\end{lem}

\begin{lem}[\cite{w2}]\label{lem3.4}
Let $T$ be a tree on $m$ edges with $x\in V_T$. Then we have
$$
D^w(x)=2D(x)-m.
$$
\end{lem}

\begin{lem}[\cite{E1}]\label{lem4.1}
Let $T$ be a tree on $n$ vertices. Then
$$
(n-1)^2\leq W(T)\leq \frac{n^3-n}{6}.
$$
The lower bound holds with equality if and only if $T\cong S_n$, whereas the upper bound holds with equality if and only if $T\cong P_n.$
\end{lem}

\begin{lem}[\cite{E2}]\label{lem4.3}
Let $C_n$ be a cycle on $n$ vertices with $x\in V_{C_n}$. Then $K\!f(C_n)=\frac{n^3-n}{12}$ and $R(x)=\frac{n^2-1}{6}$.
\end{lem}

It is known that classical distance satisfies the triangle inequality, that is, $d(x,y)+d(y,z)\geq d(x,z)$ for any $x,y,z\in V_G$. It is interesting to see that the property is enjoyed by resistance distance as well.
\begin{lem}[\cite{w8}]\label{lem3.03}
Let $G$ be a connected graph with $x,y,z\in V_G$.  Then
$$
r_G(x,y)+r_G(y,z)\geq r_G(x, z).
$$
\end{lem}

\section{Technical lemmas}\label{sec:2}

In this section we present a few technical lemmas.

\begin{lem}\label{lem3.5}
Let $G\in \mathscr{U}_{n,l}$ and $x$ be a vertex on the unique cycle. Then
$$
R^w(x)=2R(x)-(n-l).
$$
\end{lem}
\begin{proof}
We proceed by induction on the order $n$ of $G$. If $n=l$, then $G\cong C_n$ and $d(y)=2$ for any $y\in V_G$. Thus we have
$$
R^w(x)=\sum_{y\in {V_G}}d(y)r(x,y)=2\sum_{y\in {V_G}}r(x,y)=2R(x).
$$
Now suppose that $G\in \mathscr{U}_{n,l}$ with $n\geq l+1$, then $G$ contains at least one pendant vertex. Let $u$ be a pendant vertex of $G$ and $v$ be the unique neighbor of $u$ and denote $\widetilde{G}=G-u$. It is obvious that $\widetilde{G}\in \mathscr{U}_{n-1,l}$.

In what follows, we use $d_G(x)$ (resp. $ r_G(x,y), d_G(x,y), R_G(x), D_G(x), R_G^w(x), D_G^w(x)$) instead of $d(x)$ (resp. $r(x,y), d(x,y), R(x), D(x), R^w(x), D^w(x)$)  to emphasize the dependence on the graph $G$. On the one hand, it follows from Lemma \ref{lem3.3} that
\begin{eqnarray*}
     R_G^w(x)&=&\sum_{y\in {V_G}}d_G(y)r_G(x,y)\\
     &=&\sum_{y\in {V_{\widetilde{G}}\setminus\{v\}}}d_{\widetilde{G}}(y)r_{\widetilde{G}}(x,y)+(d_{\widetilde{G}}(v)+1)r_{\widetilde{G}}(x,v)+r_{\widetilde{G}}(x,v)+1\\
     &=&\sum_{y\in V_{\widetilde{G}}}d_{\widetilde{G}}(y)r_{\widetilde{G}}(x,y)+2r_{\widetilde{G}}(x,v)+1\\
     &=&R_{\widetilde{G}}^w(x)+2r_{\widetilde{G}}(x,v)+1.
\end{eqnarray*}
On the other hand,
\begin{eqnarray*}
     R_G(x)=\sum_{y\in V_{\widetilde{G}}}r_{\widetilde{G}}(x,y)+r_{\widetilde{G}}(x,v)+1=R_{\widetilde{G}}(x)+r_{\widetilde{G}}(x,v)+1.
\end{eqnarray*}
Therefore, by induction hypothesis we have
\begin{eqnarray*}
R_G^w(x)&=&R_{\widetilde{G}}^w(x)+2r_{\widetilde{G}}(x,v)+1\\
&=&2R_{\widetilde{G}}(x)-(n-1-l)+2r_{\widetilde{G}}(x,v)+1\\
&=&2R_G(x)-(n-l),
\end{eqnarray*}
as desired.
\end{proof}
More generally, we are now in a position to give a formula for $R^w(x)$ in terms of $R(x)$ for any $x\in V_G$, where $G\in \mathscr{U}_{n,l}.$
\begin{lem}\label{lem3.6}
Let $G=(C_l;T_1,T_2,\ldots, T_l)$ be in $\mathscr{U}_{n,l}$ with $x\in V_{T_k}$. Then
$$
R^w(x)=2R(x)+2d(x,v_k)-(n-l).
$$
\end{lem}
\begin{proof}
Note that
$$
\sum_{y\in V_{T_k}}d_G(y)=\sum_{y\in V_{T_k}}d_{T_k}(y)+2=2|E_{T_k}|+2=2n_k
$$
and
$$
\sum_{y\in V_{T_k}}d_G(y)d_G(x,y)=\sum_{y\in V_{T_k}}d_{T_k}(y)d_{T_k}(x,y)+2d_G(x,v_k)=D_{T_k}^w(x)+2d_G(x,v_k).
$$
Then it follows from Lemma \ref{lem3.3} that
\begin{eqnarray*}
     R_G^w(x)&=&\sum_{y\in {V_G}}d_G(y)r_G(x,y)\\
     &=&\sum_{j=1,j\neq k}^l\sum_{y\in V_{T_j}}d_G(y)(d_{G}(x,v_k)+r_G(v_k,y))+\sum_{y\in V_{T_k}}d_G(y)d_{T_k}(x,y)\\
     &=&2(n-n_k)d_G(x,v_k)+R_G^w(v_k)-\sum_{y\in V_{T_k}}d_G(y)d_G(v_k,y)+D_{T_k}^w(x)+2d_G(x,v_k)\\
     &=&2(n-n_k+1)d_G(x,v_k)+R_G^w(v_k)-D_{T_k}^w(v_k)+D_{T_k}^w(x).
\end{eqnarray*}
Note that
\begin{eqnarray*}
     R_G(x)&=&\sum_{j=1,j\neq k}^l\sum_{y\in V_{T_j}}(d_{G}(x,v_k)+r_G(v_k,y))+\sum_{y\in V_{T_k}}d_{T_k}(x,y)\\
     &=&(n-n_k)d_G(x,v_k)+R_G(v_k)-\sum_{y\in V_{T_k}}d_{T_k}(v_k,y)+\sum_{y\in V_{T_k}}d_{T_k}(x,y)\\
     &=&(n-n_k)d_G(x,v_k)+R_G(v_k)-D_{T_k}(v_k)+D_{T_k}(x).
\end{eqnarray*}
Therefore, together with Lemmas \ref{lem3.4} and \ref{lem3.5}, we have
\begin{eqnarray*}
R_G^w(x)-2 R_G(x)=2d(x,v_k)-(n-l).
\end{eqnarray*}
This completes the proof.
\end{proof}
From Lemma \ref{lem3.6}, we can prove the following explicit formula for the expected hitting time on unicyclic graphs, which will play an important role in the proof of our main results.
\begin{lem}\label{lem3.7}
Let $G\in \mathscr{U}_{n,l}$ with $x\in V_{T_i}$ and $y\in V_{T_j}$. Then
$$
H_{xy}=nr(x,y)+R(y)-R(x)+d(y,v_j)-d(x,v_i).
$$
\end{lem}
\begin{proof}
Note that  $|E_G|=|V_G|=n$. Then it follows from Lemmas \ref{lem3.1} and \ref{lem3.6} that we have
\begin{eqnarray*}
H_{xy}&=&nr(x,y)+\frac{1}{2}(R^w(y)-R^w(x))\\
&=&nr(x,y)+R(y)-R(x)+d(y,v_j)-d(x,v_i),
\end{eqnarray*}
as desired.
\end{proof}
\begin{lem}\label{lem5.1}
Let $T$ be a tree on $n$ vertices with $v\in V_T$. Then
$$
n-1\leq (2n-1)D_T(v)-2W(T)\leq \frac{n(n-1)(4n-5)}{6}.
$$
The lower bound holds with equality if and only if $T\cong S_n$ with $v$ as its center, whereas the upper bound holds with equality if and only if $T\cong P_n$ with $v$ as one of its end vertices.
\end{lem}
\begin{proof}
Let $T$ be a tree on $n$ vertices with $v\in V_T$. For $uv\in E_T$, denote $T_u$ and $T_v$ be the component of $T-uv$ containing $u$ and $v$, respectively. Assume that $|V_{T_u}|=n_u$ and $|V_{T_v}|=n_v$. Obviously, $n_u+n_v=n.$

If $n_u\geq2$, let $T'=T-\{uw_1,uw_2,\ldots,uw_k\}+\{vw_1,vw_2,\ldots,vw_k\}$, where $w_1,w_2,\ldots,w_k$ are the neighbors of $u$ in $T_u$. Then it is routine to check that $D_T(v)-D_{T'}(v)=n_u-1$ and $W(T)-W(T')=n_v(n_u-1)$. Therefore,
$$
(2n-1)D_T(v)-2W(T)-\left((2n-1)D_{T'}(v)-2W(T')\right)=(n_u-1)(2n_u-1)>0,
$$
i.e,
$$
(2n-1)D_T(v)-2W(T)>(2n-1)D_{T'}(v)-2W(T').
$$
By iterating the argument, we can get $(2n-1)D_T(v)-2W(T)$ attains its minimum (resp. maximum) if and only if $T\cong S_n$ with $v$ as its center (resp. $T\cong P_n$ with $v$ as one of its end vertices). By a direct calculation, we obtain our assertion.
\end{proof}
\section{Proofs of Theorems \ref{th2} and \ref{th1}}\label{sec:3}\setcounter{equation}{0}
First we show Theorem~\ref{th2}, which mainly follows from Lemma \ref{lem3.7}.
\subsection{The proof of Theorem~\ref{th2}}\setcounter{equation}{0}
The proofs of (i) and (ii) in Theorem \ref{th2} follow almost directly from Lemma \ref{lem3.7}, and are rather similar to each other in nature. Here we only provide the proof for (i) and the proof of (ii) is omitted here. In view of Lemma \ref{lem3.7}, we can get
\begin{eqnarray*}
     CC(x)&=&\sum_{y\in V_G}\left[nr(x,y)+R(y)-R(x)-d(x,v_k)\right]+\sum_{i=1}^l\sum_{y\in V_{T_i}}d(y,v_i)\\
     &=&nR(x)+2K\!f(G)-nR(x)-nd(x,v_k)+\sum_{i=1}^lD_{T_i}(v_i)\\
     &=&2K\!f(G)+\sum_{i=1}^lD_{T_i}(v_i)-nd(x,v_k).
\end{eqnarray*}

This completes the proof.\qed

Now we provide the proof for Theorem~\ref{th1}, which mainly follows from Lemma \ref{lem3.6}.
\subsection{The proof of Theorem~\ref{th1}}

(i)\ In view of (\ref{eq:1.1b}), we have $K\!f^+(G)=\sum_{x\in V_G}\sum_{y\in V_G}d(y)r(x,y)=\sum_{x\in V_G}R^w(x).$ Then it follows from Lemma \ref{lem3.6} that we have
\begin{eqnarray*}
     K\!f^+(G)&=&\sum_{x\in V_G}\left[2R(x)-(n-l)\right]+2\sum_{i=1}^l\sum_{x\in V_{T_i}}d(x,v_i)\\
     &=&4K\!f(G)+2\sum_{i=1}^lD_{T_i}(v_i)-n(n-l).
\end{eqnarray*}

(ii)\ By (\ref{eq:1.1a}), we can get $K\!f^*(G)=\frac{1}{2}\sum_{x\in V_G}\sum_{y\in V_G}d(x)d(y)r(x,y)=\frac{1}{2}\sum_{x\in V_G}d(x)R^w(x).$ Then it follows from Lemmas \ref{lem3.4} and \ref{lem3.6} that we have
\begin{eqnarray*}
     K\!f^*(G)&=&\frac{1}{2}\sum_{x\in V_G}d(x)\left[2R(x)-(n-l)\right]+\sum_{i=1}^l\sum_{x\in V_{T_i}}d(x)d(x,v_i)\\
     &=&\sum_{x\in V_G}\sum_{y\in V_G}d(x)r(x,y)-n(n-l)+\sum_{i=1}^lD_{T_i}^w(v_i)\\
     &=&K\!f^+(G)-n(n-l)+\sum_{i=1}^l\left[2D_{T_i}(v_i)-(n_i-1)\right]\\
     &=&4K\!f(G)+4\sum_{i=1}^lD_{T_i}(v_i)-(2n+1)(n-l).
\end{eqnarray*}
This completes the proof.\qed
\section{Proof of Theorem \ref{th4.50}}\label{sec:4}

In order to give the proof of Theorem \ref{th4.50}, we first determine sharp upper and lower bounds on $CC(x)$ of graphs among $\mathscr{U}_{n,l}$.
\begin{thm}\label{th4.4}
Let $G\in\mathscr{U}_{n,l}$ with $x\in V_G$. Then
$$
CC(x)\leq \frac{l^3}{2}-\frac{(4n+3)l^2}{6}+\frac{n(2n^2+3n-1)}{6}
$$
with equality if and only if $G\cong P_n^l$ and $x$ lies on the unique cycle $C_l$.
\end{thm}
\begin{proof}
Choose $G$ in $\mathscr{U}_{n,l}$ such that $CC(x)$ is as large as possible with $x\in V_G.$ By Theorem~\ref{th2}(i),  $x$ must lie on the unique cycle $C_l=v_1v_2\ldots v_lv_1$. Consequently,
\begin{eqnarray}\label{eq:4.1a}
CC(x)=2K\!f(G)+\sum_{i=1}^lD_{T_i}(v_i).
\end{eqnarray}
Note that
\begin{eqnarray}\label{eq:4.1}
     K\!f(G)&=&\sum_{i=1}^lW(T_i)+\sum_{1\leq i<j\leq l}\sum_{x\in V_{T_i}}\sum_{y\in V_{T_j}}r(x,y)\notag\\
     &=&\sum_{i=1}^lW(T_i)+\sum_{1\leq i<j\leq l}\sum_{x\in V_{T_i}}\sum_{y\in V_{T_j}}\left(d(x,v_i)+r(v_i,v_j)+d(y,v_j)\right)\notag\\
     &=&\sum_{i=1}^lW(T_i)+\sum_{1\leq i<j\leq l}\sum_{x\in V_{T_i}}\left(n_jd(x,v_i)+n_jr(v_i,v_j)+D_{T_j}(v_j)\right)\notag\\
     &=&\sum_{i=1}^lW(T_i)+\sum_{1\leq i<j\leq l}\left(n_jD_{T_i}(v_i)+n_in_jr(v_i,v_j)+n_iD_{T_j}(v_j)\right)\notag\\
     &=&\sum_{i=1}^lW(T_i)+\sum_{1\leq i<j\leq l}\left(n_jD_{T_i}(v_i)+n_iD_{T_j}(v_j)\right)+\frac{1}{2}\sum_{i=1}^l\sum_{j=1}^ln_in_jr(v_i,v_j).
\end{eqnarray}
Then in view of (\ref{eq:4.1a})-(\ref{eq:4.1}) we have
\begin{align}\label{eq:4.2}
     CC(x)=&2\sum_{i=1}^lW(T_i)+2\sum_{1\leq i<j\leq l}\left(n_jD_{T_i}(v_i)+n_iD_{T_j}(v_j)\right)+\sum_{i=1}^lD_{T_i}(v_i)\notag\\
     +&\sum_{i=1}^l\sum_{j=1}^ln_in_jr(v_i,v_j).
\end{align}
Let $G=U(C_l;T_1,\ldots,T_l)$ with the unique cycle $C_l=v_1v_2\ldots v_lv_1.$ Then in order to complete the proof, it suffices to show the following two claims.
\begin{claim}
$T_i\cong P_{n_i}$ and $v_i$ is one of the end vertices of $T_i$, $i=1,2, \ldots, l$.
\end{claim}
\noindent{\bf Proof of Claim 1.}\ \
By Lemma \ref{lem4.1} we know that $W(T_i)$ attains its maximum if and only if $T_i$ is a path and it is routine to check that $D_{T_i}(v_i)$ attains its  maximum if and only if $T_i$ is a path with $v_i$ being one of its end vertices. Then in view of (\ref{eq:4.2}), Claim~1 holds.
\qed

\begin{claim}
  If $n\geq l+1$, then only one member in $\{T_1,T_2, \ldots, T_l\}$ is non-trivial.
\end{claim}
\noindent{\bf Proof of Claim 2.}\ \
Suppose to the contrary that there exist two non-trivial $T_i,T_j$. By Claim~1, $T_i\cong P_{n_i}$ and $T_j\cong P_{n_j}$. Let $a$ (resp. $b$) be another end vertex different from $v_i$ (resp. $v_j$) of $T_i$ (resp. $T_j$). Assume, without loss of generality, that $R_{G}(a)+\frac{n_i}{2}\geq R_{G}(b)+\frac{n_j}{2}$. Let $c$ be the unique neighbor of $b$ and put $G'=G-cb+ab$. Then it is obvious that $G'\in \mathscr{U}_{n,l}$.

Note that
\begin{eqnarray*}
R_{G'}(b)&=&R_{G'}(a)+n-2\\
&=&R_{G}(a)+1-r_{G}(a,b)+n-2\\
&\geq& R_{G}(a)+n-1-d_{G}(a,b)>R_{G}(a).
\end{eqnarray*}
Then in view of (\ref{eq:4.1a}), we have
\begin{eqnarray*}
CC_{G'}(x)-CC_{G}(x)&=&2\left(R_{G'}(b)-R_{G}(b)\right)+n_i-(n_j-1)\\
&>&2\left(R_{G}(a)-R_{G}(b)+\frac{n_i-n_j+1}{2}\right)>0,
\end{eqnarray*}
a contradiction to the choice of $G$, which implies that all but one of  $T_i$'s are trivial.
\qed

By Claims 1-2, we have $G\cong P_n^l$. By direct calculation, we have
\begin{eqnarray*}
CC_{G}(x)&=&2\left[\frac{l^3-l}{12}+\frac{(n-l)^3-(n-l)}{6}+\sum_{k=1}^{n-l}\left(\frac{l^2-1}{6}+kl\right)\right]+\frac{(n-l)(n-l+1)}{2}\\
&=&\frac{l^3}{2}-\frac{(4n+3)l^2}{6}+\frac{n(2n^2+3n-1)}{6}.
\end{eqnarray*}
This completes the proof.
\end{proof}

Now we consider the sharp lower bound on $CC(x)$ of graphs among $\mathscr{U}_{n,l}$ in what follows.%The sharp lower bound of $CC_G(x)$ for $G\in \mathscr{U}_{n,l}$ (resp. $G\in \mathscr{U}_n$)  and the corresponding extremal graph are in the following.

\begin{thm}\label{th4.6}
Given a graph $G$ in $\mathscr{U}_{n,l}$ with $l\not=n, n\geq6,$ let $x$ be a vertex of $G$. Then
$$
CC(x)\geq -\frac{l^3}{6}+\frac{nl^2}{3}+\frac{(7-12n)l}{6}+\frac{n(6n-7)}{3}
$$
with equality if and only if $G\cong S_n^l$ and $x$ is a pendant vertex of $G$.
\end{thm}

\begin{proof}
Choose a graph $G=(C_l; T_1,T_2,\ldots,T_l)$ in $\mathscr{U}_{n,l}$ with $x\in V_G$ such that $CC(x)$ is as small as possible. Assume that $x\in V_{T_k}$, it follows from Theorem~\ref{th2}(i) that $d(x,v_k)$ attains its maximum. So we can assume that $d(x,v_k)=\varepsilon_{T_k}(v_k)$ in this case. Therefore,
\begin{eqnarray}\label{eq:4.4s}
CC(x)=2K\!f(G)+\sum_{i=1}^lD_{T_i}(v_i)-n\varepsilon_{T_k}(v_k).
\end{eqnarray}
Assume, without loss of generality, that $x\in V_{T_1}$. Then in view of (\ref{eq:4.1}) and (\ref{eq:4.4s}) we have
\begin{align}\label{eq:4.4}
CC(x)=&\ 2\sum_{j=1}^lW(T_j)+2\sum_{1\leq i<j\leq l}\left(n_jD_{T_i}(v_i)+n_iD_{T_j}(v_j)\right)+\sum_{j=1}^lD_{T_j}(v_j)\notag\\
&+\sum_{i=1}^l\sum_{j=1}^ln_in_jr(v_i,v_j)-n\varepsilon_{T_1}(v_1)\notag\\
=&\ 2W(T_1)+2\sum_{j=2}^lW(T_j)+2\sum_{j=2}^l\left(n_jD_{T_1}(v_1)+n_1D_{T_j}(v_j)\right)+D_{T_1}(v_1)+\sum_{j=2}^lD_{T_j}(v_j)\notag\\
&+2\sum_{2\leq i<j\leq l}\left(n_jD_{T_i}(v_i)+n_iD_{T_j}(v_j)\right)+\sum_{i=1}^l\sum_{j=1}^ln_in_jr(v_i,v_j)-n\varepsilon_{T_1}(v_1)\notag\\
=&\ 2W(T_1)+(2n-2n_1+1)D_{T_1}(v_1)-n\varepsilon_{T_1}(v_1)+(2n_1+1)\sum_{j=2}^lD_{T_j}(v_j)\notag\\
&+2\sum_{j=2}^lW(T_j)+2\sum_{2\leq i<j\leq l}\left(n_jD_{T_i}(v_i)+n_iD_{T_j}(v_j)\right)+\sum_{i=1}^l\sum_{j=1}^ln_in_jr(v_i,v_j).
\end{align}

In order to complete the proof, it suffices to show the following claims.
\setcounter{claim}{0}
\begin{claim}
 $2W(T_1)+(2n-2n_1+1)D_{T_1}(v_1)-n\varepsilon_{T_1}(v_1)=(2n-1)n_1-3n+1$ and $T_1\cong S_{n_1}$ with $v_1$ being its central vertex.
\end{claim}
\noindent{\bf Proof of Claim 1}\ \
We proceed by induction on $n_1$. For $n_1=2$, the statement is trivial. %Assume our result holds for all trees of order less than $n_1$.
Now let $T_1$ be a tree of order $n_1\geq3$ and $v_1$ a vertex for which $2W(T_1)+(2n-2n_1+1)D_{T_1}(v_1)-n\varepsilon_{T_1}(v_1)$ is as small as possible. Next we show that there exists a pendant vertex being adjacent to $v_1$ in $T_1$.

If $d_{T_1}(v_1)=1$, let $u_1$ be the unique neighbor of $v_1$ in $T_1$. Then we have $D_{T_1}(v_1)=D_{T_1}(u_1)+n_1-2$ and $\varepsilon_{T_1}(v_1)=\varepsilon_{T_1}(u_1)+1.$ Thus
\begin{align*}
(2n-2n_1+1)[(D_{T_1}(v_1)-n\varepsilon_{T_1}(v_1))-(D_{T_1}(u_1)-n\varepsilon_{T_1}(u_1))]%=&(2n-2n_1+1)(n_1-2)-n\\
=&-2n_1^2+(2n+5)n_1-5n-2.
\end{align*}
Let $g(t)=-2t^2+(2n+5)t-5n-2$ be a real function in $t$ with $t\in [3,n-2]$. Since $g(3)=n-5>0$ and $g(n-2)=4n-20>0$, we have $g(t)>0$ for $t\in [3,n-2]$. Hence $2W(T_1)+(2n-2n_1+1)D_{T_1}(u_1)-n\varepsilon_{T_1}(u_1)<2W(T_1)+(2n-2n_1+1)D_{T_1}(v_1)-n\varepsilon_{T_1}(v_1),$ which contradicts the choice of $v_1$.

Hence, $d_{T_1}(v_1)\geq 2$. Let $N_{T_1}(v_1)=\{u_1,u_2,\ldots,u_k\}$ and $F_i$ be the connected component of $T_1-v_1u_i$ containing $u_i,\, i=1,2,\ldots, k$. If each of $\{F_1, F_2, \ldots, F_k\}$ is non-trivial, then assume without loss of generality that $F_1$ contains a vertex $u$ such that $d_{T_1}(u,v_1)=\varepsilon_{T_1}(v_1)$. Put $T_1':=T_1-\{u_kx: \, x\in N_{F_k}(u_k)\}+\{v_1x:\, x\in N_{F_k}(u_k)\}$. Then it is routine to check that $$W(T_1)-W(T_1')=(n_1-l_k-1)(l_k-1),\ \ \ \  D_{T_1}(v_1)-D_{T_1'}(v_1)=l_k-1$$ and $\varepsilon_{T_1}(v_1)=\varepsilon_{T_1'}(v_1)$, where $l_k=|V_{F_k}|$. Consequently,
\begin{align*}
&[2W(T_1')+(2n-2n_1+1)D_{T_1'}(v_1)-n\varepsilon_{T_1'}(v_1)]-\left[2W(T_1)+(2n-2n_1+1)D_{T_1}(v_1)-n\varepsilon_{T_1}(v_1)\right]\\
=&-2(n_1-l_k-1)(l_k-1)-(2n-2n_1+1)(l_k-1)\\
=&-(l_k-1)(2n-2l_k-1)<0,
\end{align*}
a contradiction to the choice of $T_1$. Therefore, $\{F_1, F_2, \ldots, F_k\}$ contains a trivial member, say $F_k.$ Thus, $u_k$ is a pendant vertex of $T_1$ being adjacent to $v_1$. Put $T_1''=T_1-u_k$. Then it is obvious that
$$
  W(T_1)%=W(T_1'')+D_{T_1}(u_k)
  =W(T_1'')+D_{T_1''}(v_1)+n_1-1,\ \ \ \ \  D_{T_1}(v_1)=D_{T_1''}(v_1)+1
$$
and $\varepsilon_{T_1}(v_1)=\varepsilon_{T_1''}(v_1)$. By induction, one has
\begin{align}\label{eq:4.04}
2W(T_1)+(2n-2n_1+1)D_{T_1}(v_1)-n\varepsilon_{T_1}(v_1)=&\ 2W(T_1'')+(2n-2n_1+3)D_{T_1''}(v_1)\notag\\
&-n\varepsilon_{T_1''}(v_1)+2n-1\notag\\
\geq&\ (2n-1)(n_1-1)-3n+1+2n-1\\
=&\ (2n-1)n_1-3n+1.\notag
\end{align}
The equality in (\ref{eq:4.04}) holds if and only if $T_1''\cong S_{n_1-1}$ and $v_1$ is its central vertex. Therefore, $T_1\cong S_{n_1}$ with $v_1$ being its central vertex, as desired.
\qed

%Now suppose that $G=U(C_l,T_1,\ldots,T_l)$ has minimal cover cost among $\mathscr{U}_{n,l}.$
\begin{claim}
$T_i\cong S_{n_i}$ and $v_i$ is the central of $T_i$ for $2\leq i\leq l$.
\end{claim}
\noindent{\bf Proof of Claim 2}\ \
For $2\leq i \leq l$, it follows from Lemma \ref{lem4.1} that we know that $W(T_i)$ is minimal if and only if $T_i$ is a star and it is routine to check that $D_{T_i}(v_i)$ attains its minimal if and only if $T_i$ is a star with $v_i$ as its central. Then in view of (\ref{eq:4.4}), $T_i\cong S_{n_i}$ for
$2\leq i \leq l.$
\qed
\begin{claim}
All but one of $T_i$'s are trivial.
\end{claim}
\noindent{\bf Proof of Claim 3}\ \
By Claims 1-2, we have $T_i\cong S_{n_i}$ and $v_i$ is the central of $T_i$ for $1\leq i\leq l$. Then by (\ref{eq:4.4s}) we have
\begin{eqnarray}\label{eq:4.5}
CC(x)&=&2K\!f(G)+\sum_{i=1}^lD_{T_i}(v_i)-n\varepsilon_{T_i}(v_i)\notag\\
&=&2K\!f(G)+\sum_{i=1}^l(n_i-1)-n\notag\\
&=&2K\!f(G)-l.
\end{eqnarray}
Suppose to the contrary that there exist two non-trivial, say $T_i$ and $T_j$. Let $u$ and $v$ be two leaves of $T_i$ and $T_j$, respectively. Without loss of generality, assume that $R_{G}(u)\leq R_{G}(v)$. Let $G'=G-v_jv+v_iv$. Then it is obvious that $G'\in \mathscr{U}_{n,l}$.

Note that
$$
R_{G'}(v)=R_{G'}(u)=R_{G}(u)-r_{G}(u,v)+2=R_{G}(u)-r_{G}(v_i,v_j)<R_{G}(u).
$$
Then in view of (\ref{eq:4.5}) we have
$$
CC_{G'}(x)-CC_{G}(x)=2\left(K\!f(G')-K\!f(G)\right)=2\left(R_{G'}(v)-R_{G}(v)\right)<2\left(R_{G}(u)-R_{G}(v)\right)\leq0,
$$
a contradiction to the choice of $G$. Hence, our result holds.
\qed

By Claims 1-3, we have $G\cong S_n^l$ and $x$ is a pendant vertex of $G$. By direct calculation, we have
\begin{eqnarray*}
CC_{G}(x)&=&2\left[\frac{l^3-l}{12}+2\binom{n-l}{2}+(n-l)\left(\frac{l^2-1}{6}+l\right)\right]-l\\
&=&-\frac{l^3}{6}+\frac{nl^2}{3}+\frac{(7-12n)l}{6}+\frac{n(6n-7)}{3}.
\end{eqnarray*}
This completes the proof.
\end{proof}

\begin{cor}\label{lem4.7}
Given a graph $G$ in $\mathscr{U}_n$ with $l\not=n, n\geq6,$ let $x$ be a vertex of $G$. Then
%Let $G\in\mathscr{U}_n\setminus \{C_n\}$ with $n\geq6$ and $x\in V_G$. Then
\begin{eqnarray*}
CC(x)\geq\left\{
           \begin{array}{lll}
            \frac{n^3-n^2+4n-6}{6}, & \hbox{if $6\leq n\leq 8;$} \\
             2n^2-5n-6, & \hbox{if $9\leq n\leq 15;$}\\
             2n^2-\frac{16}{3}n-1, & \hbox{if $n\geq16.$}
             \end{array}
           \right.
\end{eqnarray*}
The first equality holds if and only if $G\cong S_n^{n-1}$ with $x$ being a pendant vertex of $G;$ the second equality holds if and only if $G\cong S_n^4$ with $x$ being a pendant vertex of $G$ and the last equality holds if and only if $G\cong S_n^3$ with $x$ being a pendant vertex of $G.$
\end{cor}
\begin{proof}
Let
$$
h(t)=-\frac{t^3}{6}+\frac{n}{3}t^2+\frac{7-12n}{6}t+\frac{n(6n-7)}{3}
$$
be a real function in $t$ with $t\in[3,n-1]$. Then $h'(t)=-\frac{1}{2}t^2+\frac{2n}{3}t+\frac{7-12n}{6}$. Put $\Sigma:=\frac{4n^2-36n+21}{9}$. Then we proceed by distinguishing the following two cases to show our result.

{\bf Case 1.}\ $6\leq n\leq8.$ In this case, $\Sigma<0$ and then $h'(t)<0$ for $t\in[3,n-1]$. Thus $h(t)$ is monotone decreasing on $[3,n-1]$. Hence, we have
$$
h(t)\geq h(n-1)=\frac{n^3-n^2+4n-6}{6}.
$$
Therefore, by Theorem \ref{th4.6}, our result holds in this case.

{\bf Case 2.}\ $n\geq9.$ In this case, $\Sigma>0$ and the roots of $h'(t)=0$ are $t_1=\frac{2n-\sqrt{4n^2-36n+21}}{3}$ and $t_2=\frac{2n+\sqrt{4n^2-36n+21}}{3}.$

If $n=9,$ then $4<t_1=\frac{18-\sqrt{21}}{3}<5$ and $7<t_2=\frac{18+\sqrt{21}}{3}<8$. It is routine to check that
$h'(t)<0$ for $t\in[3,t_1)\cup (t_2,8]$, and $h'(t)>0$ for $t\in(t_1,t_2)$, which implies that $h(t)$ is monotone decreasing for $t\in [3,t_1)\cup(t_2,8]$, and it is monotone increasing for $t\in (t_1,t_2)$. Combing with Theorem \ref{th4.6}, we have
$$
CC(x)\geq\min\{h(4),h(5),h(8)\}=\min\{111,111,113\}=111=h(4).
$$

If $n=10,$ then $4<t_1<5$ and $t_2>9$. It is routine to check that
$h'(t)<0$ for $t\in[3,t_1)$ and $h'(t)>0$ for $t\in(t_1,9]$, which implies that $h(t)$ is  monotone decreasing on $[3,t_1)$ and it is monotone increasing on $(t_1,9]$. Combing with Theorem \ref{th4.6}, we have
$$
CC(x)\geq\min\{h(4),h(5)\}=\min\{144,145\}=144=h(4).
$$

If $n\geq11,$ then $3<t_1<4$ and $t_2>n-1$. It is routine to check that
$h'(t)<0$ for $t\in[3,t_1)$ and $h'(t)>0$ for $t\in(t_1,n-1]$, which implies that $h(t)$ is monotone decreasing on $[3,t_1)$ and it is monotone increasing on $(t_1,n-1]$. Combing with Theorem \ref{th4.6}, we have
$$
CC(x)\geq\min\{h(3),h(4)\}=\left\{
           \begin{array}{lll}
             2n^2-5n-6=h(4), & \hbox{if $11\leq n\leq15;$}\\
             2n^2-\frac{16}{3}n-1=h(3), & \hbox{if $n\geq16.$}
             \end{array}
           \right.
$$
This completes the proof.
\end{proof}
Now we come back to the proof of Theorem \ref{th4.50}.

\noindent{\bf  Proof of Theorem \ref{th4.50}}\ \ We first determine the sharp upper bound on $CC(x)$ of $n$-vertex unicyclic graphs.
Let
$$
f(t)= \frac{t^3}{2}-\frac{(4n+3)t^2}{6}+\frac{n(2n^2+3n-1)}{6}
$$
be a real function with $t\in[3,n]$.
Then $f'(t)=\frac{3}{2}t^2-\frac{4n+3}{3}t$ and the roots of $f'(t)=0$ are $t_1=0$ and $t_2=\frac{8n+6}{9}$. Then we proceed by distinguishing the following two cases to show our result.

{\bf Case 1}.\ $n\leq6.$ In this case, $t_2\geq n$. Then $f'(t)\leq0$ and $f(t)$ is monotone decreasing on $[3,n]$. Therefore, we have
$$
f(t)\leq f(3)=\frac{2n^3+3n^2-37n+54}{6}.
$$

{\bf Case 2}.\ $n>7.$ In this case, $t_2< n$. Then $f'(t)<0$ for $t\in[3,t_2)$ and $f'(t)>0$ for $t\in(t_2,n]$, which implies that $f(t)$ is monotone decreasing on $[3,t_2)$ and it is monotone increasing on $(t_2,n]$. Since $f(3)-f(n)=\frac{2n^3+3n^2-37n+54}{6}-\frac{n^3-n}{6}=\frac{n^3+3n^2-36n+54}{6}>0$ for $n>7$, we have $f(t)\leq f(3)=\frac{2n^3+3n^2-37n+54}{6}$ in this case.

Therefore, by Cases 1-2 and Theorem \ref{th4.4} the sharp upper bound in Theorem \ref{th4.50} follows immediately.

Now we consider the sharp lower bound in Theorem \ref{th4.50}.

By Theorem \ref{th2} and Lemma \ref{lem4.3}, we have $CC_{C_n}(x)=2K\!f(C_n)=\frac{n^3-n}{6}$ with $x\in V_{C_n}.$ Let $H_i$ be the graphs with $x_i\in V_{H_i} (1\leq i\leq 5)$ as depicted in Fig. 2. Then it is routine to check that $\mathscr{U}_3=\{C_3\},\, \mathscr{U}_4=\{C_4, H_1\}$ and $\mathscr{U}_5=\{C_5, H_2, H_3, H_4, H_5\}$ and $x_i$ is the vertex of $H_i$ such that $CC_{H_i}(x_i)$ attains its minimum. By direct calculation, we have
$$
\begin{array}{cccc}
  CC_{C_3}(x)=4, & CC_{C_4}(x)=10, & CC_{H_1}(x_1)=\frac{29}{3}, & CC_{H_2}(x_2)=CC_{H_3}(x_3)=\frac{67}{3}, \\
  CC_{H_4}(x_4)=\frac{71}{3}, & CC_{H_5}(x_5)=19, & CC_{C_5}(x)=20. &
\end{array}
$$
Then the sharp lower bound in Theorem \ref{th4.50} follows immediately from Corollary~\ref{lem4.7}.
\begin{figure}[h!]
\begin{center}
\psfrag{a}{$H_1$}\psfrag{b}{$H_2$}\psfrag{c}{$H_3$}
\psfrag{d}{$H_4$}\psfrag{e}{$H_5$}
\psfrag{f}{$x_1$}\psfrag{g}{$x_2$}\psfrag{h}{$x_3$}
\psfrag{i}{$x_4$}\psfrag{j}{$x_5$}
\includegraphics[width=120mm]{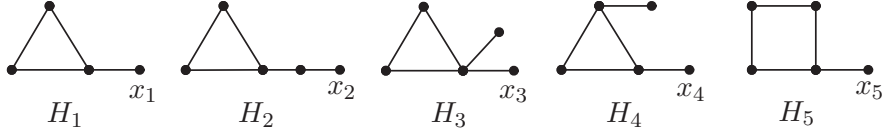} \\
\caption{The graphs $H_1, H_2, H_3, H_4$ and $H_5$.} %\label{*}
\end{center}
\end{figure}
\qed
\section{Proof of Theorem \ref{thm5}}\label{sec:5}
In order to give the proof of Theorem \ref{thm5}, we first determine sharp upper and lower bounds on $RC(x)$ of graphs among $\mathscr{U}_{n,l}$.
\begin{thm}\label{th4.9}
Let $G\in\mathscr{U}_{n,l}$ with $x\in V_G$. Then
$$
RC(x)\geq \frac{l^3}{6}-\frac{7l}{6}+n.
$$
The equality holds if and only if $G\cong S_n^l$ and $x$ is the vertex of degree $n-l+2$ in $G$.
\end{thm}
\begin{proof}
Choose $G=U(C_l;T_1,\ldots,T_l)$ in $\mathscr{U}_{n,l}$ such that $RC(x)$ is as small as possible, where $x\in V_G$. We assume, without loss of generality, that $x\in V_{T_1}$. Note that
$$
\sum_{u\in V_{T_j}}r(u,x)=\sum_{u\in V_{T_j}}(d(u,v_j)+r(v_1,v_j)+d(v_1,x))=D_{T_j}(v_j)+n_j\left(d(v_1,x)+r(v_1,v_j)\right)
$$
for $j\neq 1$.
Then it follows from Theorem \ref{th2}(ii) and (\ref{eq:4.1}) that
\begin{align}\label{eq:4.3t}
RC(x)=&2nR(x)-2K\!f(G)-\sum_{i=1}^lD_{T_i}(v_i)+nd(x, v_1)\notag\\
=&2n\left[\sum_{j=2}^l(D_{T_j}(v_j)+n_j(d(v_1,x)+r(v_1,v_j)))+D_{T_1}(x)\right]\notag\\
&-\sum_{i=1}^lD_{T_i}(v_i)+nd(x, v_1)-2K\!f(G)\notag\\
=&n(2n-2n_1+1)d(x,v_1)-D_{T_1}(v_1)+2nD_{T_1}(x)+(2n-1)\sum_{j=2}^lD_{T_j}(v_j)+2n\sum_{j=2}^lr(v_1,v_j)\notag\\
&-2\left[\sum_{i=1}^lW(T_i)+\sum_{1\leq i<j\leq l}\left(n_jD_{T_i}(v_i)+n_iD_{T_j}(v_j)\right)+\frac{1}{2}\sum_{i=1}^l\sum_{j=1}^ln_in_jr(v_i,v_j) \right]\notag\\
=&n(2n-2n_1+1)d(x,v_1)+2nD_{T_1}(x)-2W(T_1)-(2n-2n_1+1)D_{T_1}(v_1)\notag\\
&+\sum_{j=2}^l\left[(2n_j-1)D_{T_j}(v_j)-2W(T_j)\right]+2n\sum_{j=2}^lr(v_1,v_j)-\sum_{i=1}^l\sum_{j=1}^ln_in_jr(v_i,v_j).
\end{align}
For convenience, put
$$
F(T_1,v_1,x):=n(2n-2n_1+1)d(x,v_1)+2nD_{T_1}(x)-2W(T_1)-(2n-2n_1+1)D_{T_1}(v_1).
$$

By a similar discussion as in the proof of Claim 1 in Theorem \ref{th4.6}, we have the following claim.\setcounter{claim}{0}
\begin{claim}
$F(T_1,v_1,x)= n_1-1$ and $T_1\cong S_{n_1}$ with $x=v_1$ being its central vertex.
\end{claim}
%\begin{proof}
%The prove is rather similar as the Claim 1 of Theorem \ref{th4.6}, which is omitted here.
%\end{proof}

%Now suppose that $G=U(C_l,T_1,\ldots,T_l)$ has minimal reverse cover cost among $\mathscr{U}_{n,l}.$ Then
The following claim follows directly from (\ref{eq:4.3t}) and Lemma \ref{lem5.1}.
\begin{claim}
 $T_j\cong S_{n_j}$ with $v_j$ being its central vertex for $2\leq j\leq l$.
\end{claim}

By Claims 1-2, we have $T_i\cong S_{n_i}$ and $v_i$ is the central of $T_i$ for $1\leq i\leq l$. Then by Theorem~\ref{th2}(ii), we have
\begin{eqnarray}\label{eq:4.5t}
RC_{G}(x)&=&2nR_{G}(x)-2K\!f(G)-\sum_{i=1}^lD_{T_i}(v_i)+nd_{G}(v_1,x)\notag\\
&=&2nR_{G}(v_1)-2K\!f(G)-(n-l).
\end{eqnarray}
\begin{claim}
For $2\leq i\leq l, T_i$ is trivial.
\end{claim}
\noindent{\bf Proof of Claim 3}\ \
If there exist two non-trivial subtrees, say $T_i$ and $T_j.$ Let $a$ and $b$ be two leaves of $T_i$ and $T_j$, respectively. Without loss of generality, assume that $(n-1)r_{G}(v_1,v_i)-R_{G}(a)\leq (n-1)r_{G}(v_1,v_j)-R_{G}(b)$. Let $G'=G-v_jb+v_ib$. Then it is obvious that $G'\in \mathscr{U}_{n,l}$ and
\begin{eqnarray*}
R_{G'}(v_1)-R_{G}(v_1)=r_{G'}(v_1,b)-r_{G}(v_1,b)=r_{G}(v_1,v_i)-r_{G}(v_1,v_j),
\end{eqnarray*}
whereas
\begin{eqnarray*}
K\!f(G')-K\!f(G)=R_{G'}(b)-R_{G}(b)=R_{G}(a)-R_{G}(b)-r_{G}(v_i,v_j).
\end{eqnarray*}
Together with (\ref{eq:4.5t}) and Lemma \ref{lem3.03}, we have
\begin{eqnarray*}
RC_{G'}(x)-RC_{G}(x)&=&2n\left(R_{G'}(v_1)-R_{G}(v_1)\right)-2\left(K\!f(G')-K\!f(G)\right)\\
&=&2n\left(r_{G}(v_1,v_i)-r_{G}(v_1,v_j)\right)-2\left(R_{G}(a)-R_{G}(b)-r_{G}(v_i,v_j)\right)\\
&\leq&(2n+2)r_{G}(v_1,v_i)-(2n-2)r_{G}(v_1,v_j)-2R_{G}(a)+2R_{G}(b)\\
&<&2\left[(n-1)r_{G}(v_1,v_i)-R_{G}(a)-(n-1)r_{G}(v_1,v_j)+R_{G}(b)\right]\leq0,
\end{eqnarray*}
a contradiction to the choice of $G$, which implies that at most one member, say $T_k$, in $\{T_2, T_3,\linebreak \ldots, T_l\}$ is non-trivial.

Let $c,d$ be two leaves of $T_1$ and $T_k$, respectively. Let $G''=G-v_kd+v_1d$. On one hand, by a similar discussion as above, we obtain
\begin{eqnarray}\label{eq:4.01t}
RC_{G''}(x)-RC_{G}(x)=2(1-n)r_{G}(v_1,v_k)-2\left(R_{G}(c)-R_{G}(d)\right).
\end{eqnarray}
On the other hand, it is routine to check that $R_{G}(c)=R_{C_l}(v_1)+l+2(n_1-2)+(n_k-1)(2+r_{G}(v_1,v_k))$ and $R_{G}(d)=R_{C_l}(v_k)+l+2(n_k-2)+(n_1-1)(2+r_{G}(v_1,v_k))$.

Combing with (\ref{eq:4.01t}) and Lemma \ref{lem4.3}, we have
\begin{eqnarray*}
RC_{G''}(x)-RC_{G}(x)&=&2(1-n)r_{G}(v_1,v_k)-2(n_k-n_1)r_{G}(v_1,v_k)\\
&=&2r_{G}(v_1,v_k)(1+n_1-n-n_k)<0,
\end{eqnarray*}
which contradicts the choice of $G$. Therefore, all of $T_i$ are trivial for $2\leq i\leq l.$
\qed

By Claims 1-3, we have $G\cong S_n^l$ and $x$ is a vertex of degree $n-l+2$. By direct calculation, we have
\begin{eqnarray*}
RC_{G}(x)&=&2n\left(\frac{l^2-1}{6}+n-l\right)-2\left[\frac{l^3-l}{12}+2\binom{n-l}{2}+(n-l)\left(\frac{l^2-1}{6}+l\right)\right]-(n-l)\\
&=&\frac{l^3}{6}-\frac{7l}{6}+n
\end{eqnarray*}
and the proof is complete.
\end{proof}

In order to complete the proof of Theorem \ref{thm5}, we must determine the sharp upper bound on $RC(x)$ of $G$ in $\mathscr{U}_{n,l},$ where $x\in V_G.$
\begin{thm}\label{th5.4}
Let $G\in\mathscr{U}_{n,l}$ with $x\in V_G$. Then
$$
RC(x)\leq -\frac{l^3}{2}+\frac{l^2}{2}+\frac{n(n-1)(4n+1)}{6}.
$$
The equality holds if and only if $G\cong P_n^l$ and $x$ is the  vertex such that $d(x,V_{C_l})=n-l$.
\end{thm}
\begin{proof}
Choose $G=U(C_l; T_1,\ldots,T_l)$ among $\mathscr{U}_{n,l}$ such that its reverse cover cost is as small as possible. By the same notation as in the proof of Theorem 5.1, we have
$$
F(T_1,v_1,x)=n(2n-2n_1+1)d(x,v_1)+2nD_{T_1}(x)-2W(T_1)-(2n-2n_1+1)D_{T_1}(v_1).
$$
Then by a similar discussion as that of Claim 1 in Theorem \ref{th4.6}, we have \setcounter{claim}{0}
\begin{claim}
$F(T_1,v_1,x)=\frac{(n_1-1)\left[4n_1^2-(12n+5)n_1+6n(2n+1)\right]}{6}$ with $T_1\cong P_{n_1}$ and $x,\, v_1$ are two end-vertices of $T_1$.
\end{claim}
%\begin{proof}
%The prove is rather similar as the Claim 1 of Theorem \ref{th4.6}, which is omitted here.
%\end{proof}

By (\ref{eq:4.3t}) and Lemma \ref{lem5.1} the following claim follows directly.
\begin{claim}
$T_j\cong P_{n_j}$ with $v_j$ as one of its end vertices for $2\leq j\leq l$.
\end{claim}

Then by Theorem \ref{th2} (ii) and Claims 1-2, we have
\begin{eqnarray}\label{eq:5.01t}
RC_{G}(x)&=&2nR_{G}(x)-2K\!f(G)-\sum_{i=1}^lD_{T_i}(v_i)+nd_{G}(x,v_1)\notag\\
&=&2nR_{G}(x)-2K\!f(G)-\sum_{i=1}^l\frac{n_i(n_i-1)}{2}+n(n_1-1).
\end{eqnarray}
\begin{claim}
 For $2\leq i\leq l, T_i$ is trivial.
\end{claim}
\noindent{\bf Proof of Claim 3}\ \
If there exist two non-trivial members, say $T_i$ and $T_j,$ in $\{T_1, T_2, \ldots, T_l\},$ then let $u_1$ (resp. $u_2$) be another end-vertex different from $v_i$ (resp. $v_j$) of $T_i$ (resp. $T_j$). Without loss of generality, assume that $nn_i+(n+1)r_{G}(v_1,v_i)-2R_{G}(u_1)\leq nn_j+(n+1)r_{G}(v_1,v_j)-2R_{G}(u_2)$. Let $G'=G-u_2u_3+u_1u_2$, where $u_3$ is the unique neighbor of $u_2$. Then it is obvious that $G'\in \mathscr{U}_{n,l}$ and
\begin{eqnarray*}
R_{G'}(x)-R_{G}(x)=r_{G'}(x,u_2)-r_{G}(x,u_2)=r_{G}(v_1,v_i)-r_{G}(v_1,v_j)n_i-n_j+1,
\end{eqnarray*}
whereas
\begin{eqnarray*}
K\!f(G')-K\!f(G)&=&R_{G'}(u_2)-R_{G}(u_2)\\
&=&R_{G'}(u_1)-R_{G}(u_2)+n-2\\
&=&R_{G}(u_1)-R_{G}(u_2)+n-1-r_{G}(u_1,u_2)\\
&=&R_{G}(u_1)-R_{G}(u_2)-r_{G}(v_i,v_j)+n-n_i-n_j+1.
\end{eqnarray*}
Together with (\ref{eq:5.01t}) and Lemma \ref{lem3.03}, we have
\begin{eqnarray*}
RC_{G'}(x)-RC_{G}(x)&=&(2n+1)n_i-(2n-3)n_j-3-2\left(R_{G}(u_1)-R_{G}(u_2)\right)\\
&&+2n\left(r_{G}(v_1,v_i)-r_{G}(v_1,v_j)\right)+2r_{G}(v_i,v_j)\\
&\geq&(2n+1)n_i-(2n-3)n_j-3-2\left(R_{G}(u_1)-R_{G}(u_2)\right)\\
&&+(2n+2)\left(r_{G}(v_1,v_i)-r_{G}(v_1,v_j)\right)\\
&>&2\left[n(n_i-n_j)-\left(R_{G}(u_1)-R_{G}(u_2)\right)+(n+1)\left(r_{G}(v_1,v_i)-r_{G}(v_1,v_j)\right)\right]\\
&\geq&0,
\end{eqnarray*}
a contradiction, which implies that at most one member, say $T_k$, in $\{T_2, T_3, \ldots, T_l\}$ is nontrivial.

Let $u_4$ be another end-vertex different from $v_k$ of $T_k$ and $u_5$ be the unique neighbor of $u_4$.  Let $G''=G-u_4u_5+xu_4$. Then in view of (\ref{eq:5.01t}), we have \begin{eqnarray*}
  RC_{G}(x) &=& 2nR_{G}(x)-2K\!f(G)-\frac{n_1(n_1-1)+n_k(n_k-1)}{2}+n(n_1-1), \\
  RC_{G''}(u_4)&=& 2nR_{G''}(u_4)-2K\!f(G'')-\frac{n_1(n_1+1)+(n_k-1)(n_k-2)}{2}+nn_1.
\end{eqnarray*}

By a similar discussion as above, we obtain
\begin{eqnarray}\label{eq:5.02t}
R_{G''}(u_4)-R_{G}(x)=n-n_1-n_k+1-r_{G}(v_1,v_k)=l-1-r_{G}(v_1,v_k).
\end{eqnarray}
As well, it is routine to check that
\begin{eqnarray}\label{eq:5.03t}
K\!f(G'')-K\!f(G)=R_{G''}(u_4)-R_{G}(u_4)=(n-n_k+1)(l-1-r_{G}(v_1,v_k)).
\end{eqnarray}
Combing with (\ref{eq:5.02t})-(\ref{eq:5.03t}), we have
\begin{eqnarray*}
RC_{G''}(u_4)-RC_{G}(x)=(n-n_1+n_k-1)(2l-2r_{G}(v_1,v_k)-1)>0,
\end{eqnarray*}
which contradicts the choice of $G$ and $x$. Therefore, all of $T_i$ are trivial for $2\leq i\leq l.$
\qed

By Claims 1-3, we have that $G\cong P_n^l$ and contains the vertex $x$ satisfying $d(x,V_{C_l})=n-l$. By direct calculation, we have
\begin{eqnarray*}
RC_{G}(x)&=&2n\left[\frac{(n-l)(n-l-1)}{2}+\frac{l^2-1}{6}+l(n-l)\right]-2\left[\frac{l^3-l}{12}+\frac{(n-l)^3-(n-l)}{6}\right.\\
&&+\left.\sum_{k=1}^{n-l}\left(\frac{l^2-1}{6}+kl\right)\right]-\frac{(n-l)(n-l+1)}{2}+n(n-l)\\
&=&-\frac{l^3}{2}+\frac{l^2}{2}+\frac{n(n-1)(4n+1)}{6},
\end{eqnarray*}
as desired.
\end{proof}

%A straightforward consequence of Theorem \ref{th4.9} is stated in the next result.
\noindent{\bf Proof of Theorem \ref{thm5}}\ \ On the one hand, in view of Theorem \ref{th4.9}, let
$$
p(t)= \frac{t^3}{6}-\frac{7t}{6}+n
$$
be a real function in $t$, where $t\in[3,n]$. Thus $p'(t)=\frac{t^2}{2}-\frac{7}{6}\geq\frac{10}{3}$,  which implies that $p(t)$ is monotone increasing on $[3,n]$. Therefore, we have
$$
p(t)\geq p(3)=n+1.
$$
Together with Theorem \ref{th4.9}, the sharp lower bound on $RC_G(x)$ holds among $\mathscr{U}_n$.

On the other hand, in view of Theorem~\ref{th5.4}, consider the real function
$$
q(k)=-\frac{k^3}{2}+\frac{k^2}{2}+\frac{n(n-1)(4n+1)}{6}
$$
in $k$ with $k\in[3,n]$. Then $q'(k)=-\frac{k}{2}(3k-2)<0$,  which implies that $q(k)$ is monotone decreasing on $[3,n]$. Therefore, we have
$$
q(k)\leq q(3)=\frac{n(n-1)(4n+1)}{6}-9.
$$
Together with Theorem \ref{th5.4}, the sharp upper bound on $RC_G(x)$ holds among $\mathscr{U}_n$.

This completes the proof.
\qed

\end{document}